

\input amstex
\expandafter\ifx\csname mathdefs.tex\endcsname\relax
  \expandafter\gdef\csname mathdefs.tex\endcsname{}
\else \message{Hey!  Apparently you were trying to
  \string\input{mathdefs.tex} twice.   This does not make sense.} 
\errmessage{Please edit your file (probably \jobname.tex) and remove
any duplicate ``\string\input'' lines}\endinput\fi




\catcode`\X=12\catcode`\@=11

\def\n@wcount{\alloc@0\count\countdef\insc@unt}
\def\n@wwrite{\alloc@7\write\chardef\sixt@@n}
\def\n@wread{\alloc@6\read\chardef\sixt@@n}
\def\r@s@t{\relax}\def\v@idline{\par}\def\@mputate#1/{#1}
\def\l@c@l#1X{\firstpart.#1}\def\gl@b@l#1X{#1}\def\t@d@l#1X{{}}

\def\crossrefs#1{\ifx\all#1\let\tr@ce=\all\else\def\tr@ce{#1,}\fi
   \n@wwrite\cit@tionsout\openout\cit@tionsout=\jobname.cit 
   \write\cit@tionsout{\tr@ce}\expandafter\setfl@gs\tr@ce,}
\def\setfl@gs#1,{\def\@{#1}\ifx\@\empty\let\next=\relax
   \else\let\next=\setfl@gs\expandafter\xdef
   \csname#1tr@cetrue\endcsname{}\fi\next}
\def\m@ketag#1#2{\expandafter\n@wcount\csname#2tagno\endcsname
     \csname#2tagno\endcsname=0\let\tail=\all\xdef\all{\tail#2,}
   \ifx#1\l@c@l\let\tail=\r@s@t\xdef\r@s@t{\csname#2tagno\endcsname=0\tail}\fi
   \expandafter\gdef\csname#2cite\endcsname##1{\expandafter
     \ifx\csname#2tag##1\endcsname\relax?\else\csname#2tag##1\endcsname\fi
     \expandafter\ifx\csname#2tr@cetrue\endcsname\relax\else
     \write\cit@tionsout{#2tag ##1 cited on page \folio.}\fi}
   \expandafter\gdef\csname#2page\endcsname##1{\expandafter
     \ifx\csname#2page##1\endcsname\relax?\else\csname#2page##1\endcsname\fi
     \expandafter\ifx\csname#2tr@cetrue\endcsname\relax\else
     \write\cit@tionsout{#2tag ##1 cited on page \folio.}\fi}
   \expandafter\gdef\csname#2tag\endcsname##1{\expandafter
      \ifx\csname#2check##1\endcsname\relax
      \expandafter\xdef\csname#2check##1\endcsname{}%
      \else\immediate\write16{Warning: #2tag ##1 used more than once.}\fi
      \multit@g{#1}{#2}##1/X%
      \write\t@gsout{#2tag ##1 assigned number \csname#2tag##1\endcsname\space
      on page \number\count0.}%
   \csname#2tag##1\endcsname}}

\def\multit@g#1#2#3/#4X{\def\t@mp{#4}\ifx\t@mp\empty%
      \global\advance\csname#2tagno\endcsname by 1 
      \expandafter\xdef\csname#2tag#3\endcsname
      {#1\number\csname#2tagno\endcsnameX}%
   \else\expandafter\ifx\csname#2last#3\endcsname\relax
      \expandafter\n@wcount\csname#2last#3\endcsname
      \global\advance\csname#2tagno\endcsname by 1 
      \expandafter\xdef\csname#2tag#3\endcsname
      {#1\number\csname#2tagno\endcsnameX}
      \write\t@gsout{#2tag #3 assigned number \csname#2tag#3\endcsname\space
      on page \number\count0.}\fi
   \global\advance\csname#2last#3\endcsname by 1
   \def\t@mp{\expandafter\xdef\csname#2tag#3/}%
   \expandafter\t@mp\@mputate#4\endcsname
   {\csname#2tag#3\endcsname\lastpart{\csname#2last#3\endcsname}}\fi}
\def\t@gs#1{\def\all{}\m@ketag#1e\m@ketag#1s\m@ketag\t@d@l p
\let\realscite\scite
\let\realstag\stag
   \m@ketag\gl@b@l r \n@wread\t@gsin
   \openin\t@gsin=\jobname.tgs \re@der \closein\t@gsin
   \n@wwrite\t@gsout\openout\t@gsout=\jobname.tgs }
\outer\def\localtags{\t@gs\l@c@l}
\outer\def\globaltags{\t@gs\gl@b@l}
\outer\def\newlocaltag#1{\m@ketag\l@c@l{#1}}
\outer\def\newglobaltag#1{\m@ketag\gl@b@l{#1}}

\newif\ifpr@ 
\def\m@kecs #1tag #2 assigned number #3 on page #4.%
   {\expandafter\gdef\csname#1tag#2\endcsname{#3}
   \expandafter\gdef\csname#1page#2\endcsname{#4}
   \ifpr@\expandafter\xdef\csname#1check#2\endcsname{}\fi}
\def\re@der{\ifeof\t@gsin\let\next=\relax\else
   \read\t@gsin to\t@gline\ifx\t@gline\v@idline\else
   \expandafter\m@kecs \t@gline\fi\let \next=\re@der\fi\next}
\def\pretags#1{\pr@true\pret@gs#1,,}
\def\pret@gs#1,{\def\@{#1}\ifx\@\empty\let\n@xtfile=\relax
   \else\let\n@xtfile=\pret@gs \openin\t@gsin=#1.tgs \message{#1} \re@der 
   \closein\t@gsin\fi \n@xtfile}

\newcount\sectno\sectno=0\newcount\subsectno\subsectno=0
\newif\ifultr@local \def\ultralocal{\ultr@localtrue}
\def\firstpart{\number\sectno}
\def\lastpart#1{\ifcase#1 \or a\or b\or c\or d\or e\or f\or g\or h\or 
   i\or k\or l\or m\or n\or o\or p\or q\or r\or s\or t\or u\or v\or w\or 
   x\or y\or z \fi}

\def\resetall{\global\advance\sectno by 1\subsectno=0
   \gdef\firstpart{\number\sectno}\r@s@t}
\def\resetsub{\global\advance\subsectno by 1
   \gdef\firstpart{\number\sectno.\number\subsectno}\r@s@t}
\def\newsection#1\par{\resetall\vskip0pt plus.3\vsize\penalty-250
   \vskip0pt plus-.3\vsize\bigskip\bigskip
   \message{#1}\leftline{\bf#1}\nobreak\bigskip}
\def\subsection#1\par{\ifultr@local\resetsub\fi
   \vskip0pt plus.2\vsize\penalty-250\vskip0pt plus-.2\vsize
   \bigskip\smallskip\message{#1}\leftline{\bf#1}\nobreak\medskip}


\newdimen\marginshift

\newdimen\margindelta
\newdimen\marginmax
\newdimen\marginmin

\def\margininit{       
\marginmax=3 true cm                  
				      
\margindelta=0.1 true cm              
\marginmin=0.1true cm                 
\marginshift=\marginmin
}    

\def\t@gsjj#1,{\def\@{#1}\ifx\@\empty\let\next=\relax\else\let\next=\t@gsjj
   \def\@@{p}\ifx\@\@@\else
   \expandafter\gdef\csname#1cite\endcsname##1{\citejj{##1}}
   \expandafter\gdef\csname#1page\endcsname##1{?}
   \expandafter\gdef\csname#1tag\endcsname##1{\tagjj{##1}}\fi\fi\next}
\newif\ifshowstuffinmargin
\showstuffinmarginfalse
\def\jjtags{\showstuffinmargintrue
\ifx\all\relax\else\expandafter\t@gsjj\all,\fi}

\def\tagjj#1{\realstag{#1}\mginpar{\zeigen{#1}}}
\def\citejj#1{\zeigen{#1}\mginpar{\rechnen{#1}}}

\def\rechnen#1{\expandafter\ifx\csname stag#1\endcsname\relax ??\else
                           \csname stag#1\endcsname\fi}

\newdimen\theight

\def\marginfont{\sevenrm}

\def\trymarginbox#1{\setbox0=\hbox{\marginfont\hskip\marginshift #1}%
		\global\marginshift\wd0 
		\global\advance\marginshift\margindelta}

\def \mginpar#1{%
\ifvmode\setbox0\hbox to \hsize{\hfill\rlap{\marginfont\quad#1}}%
\ht0 0cm
\dp0 0cm
\box0\vskip-\baselineskip
\else 
             \vadjust{\trymarginbox{#1}%
		\ifdim\marginshift>\marginmax \global\marginshift\marginmin
			\trymarginbox{#1}%
                \fi
             \theight=\ht0
             \advance\theight by \dp0    \advance\theight by \lineskip
             \kern -\theight \vbox to \theight{\rightline{\rlap{\box0}}%
\vss}}\fi}


\def\t@gsoff#1,{\def\@{#1}\ifx\@\empty\let\next=\relax\else\let\next=\t@gsoff
   \def\@@{p}\ifx\@\@@\else
   \expandafter\gdef\csname#1cite\endcsname##1{\zeigen{##1}}
   \expandafter\gdef\csname#1page\endcsname##1{?}
   \expandafter\gdef\csname#1tag\endcsname##1{\zeigen{##1}}\fi\fi\next}
\def\verbatimtags{\showstuffinmarginfalse
\ifx\all\relax\else\expandafter\t@gsoff\all,\fi}
\def\zeigen#1{\hbox{$\langle$}#1\hbox{$\rangle$}}
\def\margincite#1{\ifshowstuffinmargin\mginpar{\rechnen{#1}}\fi}

\def\(#1){\edef\dot@g{\ifmmode\ifinner(\hbox{\noexpand\etag{#1}})
   \else\noexpand\eqno(\hbox{\noexpand\etag{#1}})\fi
   \else(\noexpand\ecite{#1})\fi}\dot@g}

\newif\ifbr@ck
\def\eat#1{}
\def\[#1]{\br@cktrue[\br@cket#1'X]}
\def\br@cket#1'#2X{\def\temp{#2}\ifx\temp\empty\let\next\eat
   \else\let\next\br@cket\fi
   \ifbr@ck\br@ckfalse\br@ck@t#1,X\else\br@cktrue#1\fi\next#2X}
\def\br@ck@t#1,#2X{\def\temp{#2}\ifx\temp\empty\let\neext\eat
   \else\let\neext\br@ck@t\def\temp{,}\fi
   \def\teemp{#1}\ifx\teemp\empty\else\rcite{#1}\fi\temp\neext#2X}
\def\resetbr@cket{\gdef\[##1]{[\rtag{##1}]}}
\def\references{\resetbr@cket\newsection References\par}

\newtoks\symb@ls\newtoks\s@mb@ls\newtoks\p@gelist\n@wcount\ftn@mber
    \ftn@mber=1\newif\ifftn@mbers\ftn@mbersfalse\newif\ifbyp@ge\byp@gefalse
\def\defm@rk{\ifftn@mbers\n@mberm@rk\else\symb@lm@rk\fi}
\def\n@mberm@rk{\xdef\m@rk{{\the\ftn@mber}}%
    \global\advance\ftn@mber by 1 }
\def\rot@te#1{\let\temp=#1\global#1=\expandafter\r@t@te\the\temp,X}
\def\r@t@te#1,#2X{{#2#1}\xdef\m@rk{{#1}}}
\def\b@@st#1{{$^{#1}$}}\def\str@p#1{#1}
\def\symb@lm@rk{\ifbyp@ge\rot@te\p@gelist\ifnum\expandafter\str@p\m@rk=1 
    \s@mb@ls=\symb@ls\fi\write\f@nsout{\number\count0}\fi \rot@te\s@mb@ls}
\def\byp@ge{\byp@getrue\n@wwrite\f@nsin\openin\f@nsin=\jobname.fns 
    \n@wcount\currentp@ge\currentp@ge=0\p@gelist={0}
    \re@dfns\closein\f@nsin\rot@te\p@gelist
    \n@wread\f@nsout\openout\f@nsout=\jobname.fns }
\def\m@kelist#1X#2{{#1,#2}}
\def\re@dfns{\ifeof\f@nsin\let\next=\relax\else\read\f@nsin to \f@nline
    \ifx\f@nline\v@idline\else\let\t@mplist=\p@gelist
    \ifnum\currentp@ge=\f@nline
    \global\p@gelist=\expandafter\m@kelist\the\t@mplistX0
    \else\currentp@ge=\f@nline
    \global\p@gelist=\expandafter\m@kelist\the\t@mplistX1\fi\fi
    \let\next=\re@dfns\fi\next}
\def\symbols#1{\symb@ls={#1}\s@mb@ls=\symb@ls} 
\def\bigsymbol{\textstyle}
\symbols{\bigsymbol\ast,\dagger,\ddagger,\sharp,\flat,\natural,\star}
\def\ftnumbers{\ftn@mberstrue} \def\ftsymbols{\ftn@mbersfalse}
\def\paginal{\byp@ge} \def\resetftnumbers{\ftn@mber=1}
\def\ftnote#1{\defm@rk\expandafter\expandafter\expandafter\footnote
    \expandafter\b@@st\m@rk{#1}}

\long\def\jump#1\endjump{}
\def\ssum{\mathop{\lower .1em\hbox{$\textstyle\Sigma$}}\nolimits}

\def\qed{\nobreak\kern 1em \vrule height .5em width .5em depth 0em}
\def\newneq{\hbox{\rlap{\hbox to 1\wd9{\hss$=$\hss}}\raise .1em 
   \hbox to 1\wd9{\hss$\scriptscriptstyle/$\hss}}}
\def\subsetne{\setbox9 = \hbox{$\subset$}\mathrel{\hbox{\rlap
   {\lower .4em \newneq}\raise .13em \hbox{$\subset$}}}}
\def\supsetne{\setbox9 = \hbox{$\subset$}\mathrel{\hbox{\rlap
   {\lower .4em \newneq}\raise .13em \hbox{$\supset$}}}}

\def\vbar{\mathchoice{\vrule height6.3ptdepth-.5ptwidth.8pt\kern-.8pt}
   {\vrule height6.3ptdepth-.5ptwidth.8pt\kern-.8pt}
   {\vrule height4.1ptdepth-.35ptwidth.6pt\kern-.6pt}
   {\vrule height3.1ptdepth-.25ptwidth.5pt\kern-.5pt}}
\def\f@dge{\mathchoice{}{}{\mkern.5mu}{\mkern.8mu}}
\def\b@c#1#2{{\rm \mkern#2mu\vbar\mkern-#2mu#1}}
\def\b@b#1{{\rm I\mkern-3.5mu #1}}
\def\b@a#1#2{{\rm #1\mkern-#2mu\f@dge #1}}
\def\bb#1{{\count4=`#1 \advance\count4by-64 \ifcase\count4\or\b@a A{11.5}\or
   \b@b B\or\b@c C{5}\or\b@b D\or\b@b E\or\b@b F \or\b@c G{5}\or\b@b H\or
   \b@b I\or\b@c J{3}\or\b@b K\or\b@b L \or\b@b M\or\b@b N\or\b@c O{5} \or
   \b@b P\or\b@c Q{5}\or\b@b R\or\b@a S{8}\or\b@a T{10.5}\or\b@c U{5}\or
   \b@a V{12}\or\b@a W{16.5}\or\b@a X{11}\or\b@a Y{11.7}\or\b@a Z{7.5}\fi}}

\catcode`\X=11 \catcode`\@=12


\expandafter\ifx\csname citeadd.tex\endcsname\relax
\expandafter\gdef\csname citeadd.tex\endcsname{}
\else \message{Hey!  Apparently you were trying to
\string\input{citeadd.tex} twice.   This does not make sense.} 
\errmessage{Please edit your file (probably \jobname.tex) and remove
any duplicate ``\string\input'' lines}\endinput\fi

\sectno=-1   
\localtags
\NoBlackBoxes
\define\mr{\medskip\roster}
\define\sn{\smallskip\noindent}
\define\mn{\medskip\noindent}
\define\bn{\bigskip\noindent}
\define\ub{\underbar}

\define \nl{\newline}
\magnification=\magstep 1
\documentstyle{amsppt}

{    
\catcode`@11

\ifx\alicetwothousandloaded@\relax
  \endinput\else\global\let\alicetwothousandloaded@\relax\fi

\gdef\subjclass{\let\savedef@\subjclass
 \def\subjclass##1\endsubjclass{\let\subjclass\savedef@
   \toks@{\def\usualspace{{\rm\enspace}}\eightpoint}%
   \toks@@{##1\unskip.}%
   \edef\thesubjclass@{\the\toks@
     \frills@{{\noexpand\rm2000 {\noexpand\it Mathematics Subject
       Classification}.\noexpand\enspace}}%
     \the\toks@@}}%
  \nofrillscheck\subjclass}
} 

\pageheight{8.5truein}
\topmatter
\title{A partition relation using strongly compact cardinals \\
 Sh761} \endtitle
\author {Saharon Shelah \thanks {\null\newline I would like to thank 
Alice Leonhardt for the beautiful typing. \null\newline
Research of the author was partially  supported by the 
United States-Israel Binational Science Foundation. \null\newline
Publ.761. \null\newline
Latest Revision - 01/March/20} \endthanks} \endauthor 
\affil{Institute of Mathematics\\
 The Hebrew University\\
 Jerusalem, Israel
 \medskip
 Rutgers University\\
 Mathematics Department\\ 
New Brunswick, NJ  USA} \endaffil
\medskip
\abstract
If $\kappa$ is strongly compact and $\lambda > \kappa$ is regular,
then $\left(2^{<\lambda}\right)^+ \rightarrow (\lambda + \zeta)^2_\theta$ 
holds for $\zeta,\theta < \kappa$.  \endabstract
\endtopmatter
\document  

\expandafter\ifx\csname alice2jlem.tex\endcsname\relax
  \expandafter\xdef\csname alice2jlem.tex\endcsname{\the\catcode`@}
\else \message{Hey!  Apparently you were trying to
\string\input{alice2jlem.tex}  twice.   This does not make sense.}
\errmessage{Please edit your file (probably \jobname.tex) and remove
any duplicate ``\string\input'' lines}\endinput\fi

\expandafter\ifx\csname bib4plain.tex\endcsname\relax
  \expandafter\gdef\csname bib4plain.tex\endcsname{}
\else \message{Hey!  Apparently you were trying to \string\input
  bib4plain.tex twice.   This does not make sense.}
\errmessage{Please edit your file (probably \jobname.tex) and remove
any duplicate ``\string\input'' lines}\endinput\fi

\def\renewcommand{\newcommand}	       
\edef\cite{\the\catcode`@}%
\catcode`@ = 11
\let\@oldatcatcode = \cite
\chardef\@letter = 11
\chardef\@other = 12
%
%
%
%
\def\@innerdef#1#2{\edef#1{\expandafter\noexpand\csname #2\endcsname}}%
%
%
\@innerdef\@innernewcount{newcount}%
\@innerdef\@innernewdimen{newdimen}%
\@innerdef\@innernewif{newif}%
\@innerdef\@innernewwrite{newwrite}%
%
%
%
\def\@gobble#1{}%
%
%
%
\ifx\inputlineno\@undefined
   \let\@linenumber = \empty 
\else
   \def\@linenumber{\the\inputlineno:\space}%
\fi
%
%
%
\def\@futurenonspacelet#1{\def\cs{#1}%
   \afterassignment\@stepone\let\@nexttoken=
}%
\begingroup 
\def\\{\global\let\@stoken= }%
\\ 
\endgroup
\def\@stepone{\expandafter\futurelet\cs\@steptwo}%
\def\@steptwo{\expandafter\ifx\cs\@stoken\let\@@next=\@stepthree
   \else\let\@@next=\@nexttoken\fi \@@next}%
\def\@stepthree{\afterassignment\@stepone\let\@@next= }%
%
%
%
\def\@getoptionalarg#1{%
   \let\@optionaltemp = #1%
   \let\@optionalnext = \relax
   \@futurenonspacelet\@optionalnext\@bracketcheck
}%
%
%
\def\@bracketcheck{%
   \ifx [\@optionalnext
      \expandafter\@@getoptionalarg
   \else
      \let\@optionalarg = \empty
      \expandafter\@optionaltemp
   \fi
}%
\def\@@getoptionalarg[#1]{%
   \def\@optionalarg{#1}%
   \@optionaltemp
}%
%
%
%
\def\@nnil{\@nil}%
\def\@fornoop#1\@@#2#3{}%
\def\@for#1:=#2\do#3{%
   \edef\@fortmp{#2}%
   \ifx\@fortmp\empty \else
      \expandafter\@forloop#2,\@nil,\@nil\@@#1{#3}%
   \fi
}%
\def\@forloop#1,#2,#3\@@#4#5{\def#4{#1}\ifx #4\@nnil \else
       #5\def#4{#2}\ifx #4\@nnil \else#5\@iforloop #3\@@#4{#5}\fi\fi
}%
\def\@iforloop#1,#2\@@#3#4{\def#3{#1}\ifx #3\@nnil
       \let\@nextwhile=\@fornoop \else
      #4\relax\let\@nextwhile=\@iforloop\fi\@nextwhile#2\@@#3{#4}%
}%
%
%
%
\@innernewif\if@fileexists
\def\@testfileexistence{\@getoptionalarg\@finishtestfileexistence}%
\def\@finishtestfileexistence#1{%
   \begingroup
      \def\extension{#1}%
      \immediate\openin0 =
         \ifx\@optionalarg\empty\jobname\else\@optionalarg\fi
         \ifx\extension\empty \else .#1\fi
         \space
      \ifeof 0
         \global\@fileexistsfalse
      \else
         \global\@fileexiststrue
      \fi
      \immediate\closein0
   \endgroup
}%
%
%
%
%
\def\bibliographystyle#1{%
   \@readauxfile
   \@writeaux{\string\bibstyle{#1}}%
}%
\let\bibstyle = \@gobble
%
%
\let\bblfilebasename = \jobname
\def\bibliography#1{%
   \@readauxfile
   \@writeaux{\string\bibdata{#1}}%
   \@testfileexistence[\bblfilebasename]{bbl}%
   \if@fileexists
      \nobreak
      \@readbblfile
   \fi
}%
\let\bibdata = \@gobble
%
%
\def\nocite#1{%
   \@readauxfile
   \@writeaux{\string\citation{#1}}%
}%
\@innernewif\if@notfirstcitation
%
%
\def\cite{\@getoptionalarg\@cite}%
%
%
\def\@cite#1{%
   \let\@citenotetext = \@optionalarg
   \printcitestart
   \nocite{#1}%
   \@notfirstcitationfalse
   \@for \@citation :=#1\do
   {%
      \expandafter\@onecitation\@citation\@@
   }%
   \ifx\empty\@citenotetext\else
      \printcitenote{\@citenotetext}%
   \fi
   \printcitefinish
}%
\def\@onecitation#1\@@{%
   \if@notfirstcitation
      \printbetweencitations
   \fi
   \expandafter \ifx \csname\@citelabel{#1}\endcsname \relax
      \if@citewarning
         \message{\@linenumber Undefined citation `#1'.}%
      \fi
      \expandafter\gdef\csname\@citelabel{#1}\endcsname{%
\strut
\vadjust{\vskip-\dp\strutbox
\vbox to 0pt{\vss\parindent0cm \leftskip=\hsize 
\advance\leftskip3mm
\advance\hsize 4cm\strut\openup-4pt 
\rightskip 0cm plus 1cm minus 0.5cm ?  #1 ?\strut}}
         {\tt
            \escapechar = -1
            \nobreak\hskip0pt
            \expandafter\string\csname#1\endcsname
            \nobreak\hskip0pt
         }%
      }%
   \fi
   \csname\@citelabel{#1}\endcsname
   \@notfirstcitationtrue
}%
%
%
\def\@citelabel#1{b@#1}%
%
%
\def\@citedef#1#2{\expandafter\gdef\csname\@citelabel{#1}\endcsname{#2}}%
%
%
%
\def\@readbblfile{%
   \ifx\@itemnum\@undefined
      \@innernewcount\@itemnum
   \fi
   \begingroup
      \def\begin##1##2{%
         \setbox0 = \hbox{\biblabelcontents{##2}}%
         \biblabelwidth = \wd0
      }%
      \def\end##1{}
      %
      %
      \@itemnum = 0
      \def\bibitem{\@getoptionalarg\@bibitem}%
      \def\@bibitem{%
         \ifx\@optionalarg\empty
            \expandafter\@numberedbibitem
         \else
            \expandafter\@alphabibitem
         \fi
      }%
      \def\@alphabibitem##1{%
         \expandafter \xdef\csname\@citelabel{##1}\endcsname {\@optionalarg}%
         \ifx\biblabelprecontents\@undefined
            \let\biblabelprecontents = \relax
         \fi
         \ifx\biblabelpostcontents\@undefined
            \let\biblabelpostcontents = \hss
         \fi
         \@finishbibitem{##1}%
      }%
      \def\@numberedbibitem##1{%
         \advance\@itemnum by 1
         \expandafter \xdef\csname\@citelabel{##1}\endcsname{\number\@itemnum}%
         \ifx\biblabelprecontents\@undefined
            \let\biblabelprecontents = \hss
         \fi
         \ifx\biblabelpostcontents\@undefined
            \let\biblabelpostcontents = \relax
         \fi
         \@finishbibitem{##1}%
      }%
      \def\@finishbibitem##1{%
         \biblabelprint{\csname\@citelabel{##1}\endcsname}%
         \@writeaux{\string\@citedef{##1}{\csname\@citelabel{##1}\endcsname}}%
         \ignorespaces
      }%
      %
      %
      \let\em = \bblem
      \let\newblock = \bblnewblock
      \let\sc = \bblsc
      \frenchspacing
      \clubpenalty = 4000 \widowpenalty = 4000
      \tolerance = 10000 \hfuzz = .5pt
      \everypar = {\hangindent = \biblabelwidth
                      \advance\hangindent by \biblabelextraspace}%
      \bblrm
      \parskip = 1.5ex plus .5ex minus .5ex
      \biblabelextraspace = .5em
      \bblhook
      \input \bblfilebasename.bbl
   \endgroup
}%
%
%
\@innernewdimen\biblabelwidth
\@innernewdimen\biblabelextraspace
%
%
%
\def\biblabelprint#1{%
   \noindent
   \hbox to \biblabelwidth{%
      \biblabelprecontents
      \biblabelcontents{#1}%
      \biblabelpostcontents
   }%
   \kern\biblabelextraspace
}%
%
%
%
\def\biblabelcontents#1{{\bblrm [#1]}}%
%
%
\def\bblrm{\rm}%
%
%
\def\bblem{\it}%
%
%
\def\bblsc{\ifx\@scfont\@undefined
              \font\@scfont = cmcsc10
           \fi
           \@scfont
}%
%
%
\def\bblnewblock{\hskip .11em plus .33em minus .07em }%
%
%
\let\bblhook = \empty
%
%
%
\def\printcitestart{[}
\def\printcitefinish{]}
\def\printbetweencitations{, }
\def\printcitenote#1{, #1}
%
%
%
\let\citation = \@gobble
%
%
%
\@innernewcount\@numparams
%
%
\def\newcommand#1{%
   \def\@commandname{#1}%
   \@getoptionalarg\@continuenewcommand
}%
%
%
\def\@continuenewcommand{%
   \@numparams = \ifx\@optionalarg\empty 0\else\@optionalarg \fi \relax
   \@newcommand
}%
%
%
\def\@newcommand#1{%
   \def\@startdef{\expandafter\edef\@commandname}%
   \ifnum\@numparams=0
      \let\@paramdef = \empty
   \else
      \ifnum\@numparams>9
         \errmessage{\the\@numparams\space is too many parameters}%
      \else
         \ifnum\@numparams<0
            \errmessage{\the\@numparams\space is too few parameters}%
         \else
            \edef\@paramdef{%
               \ifcase\@numparams
                  \empty  No arguments.
               \or ####1%
               \or ####1####2%
               \or ####1####2####3%
               \or ####1####2####3####4%
               \or ####1####2####3####4####5%
               \or ####1####2####3####4####5####6%
               \or ####1####2####3####4####5####6####7%
               \or ####1####2####3####4####5####6####7####8%
               \or ####1####2####3####4####5####6####7####8####9%
               \fi
            }%
         \fi
      \fi
   \fi
   \expandafter\@startdef\@paramdef{#1}%
}%
%
%
%
%
\def\@readauxfile{%
   \if@auxfiledone \else 
      \global\@auxfiledonetrue
      \@testfileexistence{aux}%
      \if@fileexists
         \begingroup
            \endlinechar = -1
            \catcode`@ = 11
            \input \jobname.aux
         \endgroup
      \else
         \message{\@undefinedmessage}%
         \global\@citewarningfalse
      \fi
      \immediate\openout\@auxfile = \jobname.aux
   \fi
}%
%
%
\newif\if@auxfiledone
\ifx\noauxfile\@undefined \else \@auxfiledonetrue\fi
%
%
%
%
\@innernewwrite\@auxfile
\def\@writeaux#1{\ifx\noauxfile\@undefined \write\@auxfile{#1}\fi}%
%
%
%
\ifx\@undefinedmessage\@undefined
   \def\@undefinedmessage{No .aux file; I won't give you warnings about
                          undefined citations.}%
\fi
%
%
\@innernewif\if@citewarning
\ifx\noauxfile\@undefined \@citewarningtrue\fi
%
%
%
\catcode`@ = \@oldatcatcode


\def\widestnumber#1#2{}

\def\rm{\fam0 \tenrm}

\def\fakesubhead#1\endsubhead{\bigskip\noindent{\bf#1}\par}



%
%
%

%

\font\textrsfs=rsfs10
\font\scriptrsfs=rsfs7
\font\scriptscriptrsfs=rsfs5

\newfam\rsfsfam
\textfont\rsfsfam=\textrsfs
\scriptfont\rsfsfam=\scriptrsfs
\scriptscriptfont\rsfsfam=\scriptscriptrsfs

\edef\oldcatcodeofat{\the\catcode`\@}
\catcode`\@11

\def\Cal@@#1{\noaccents@ \fam \rsfsfam #1}

\catcode`\@\oldcatcodeofat


\expandafter\ifx \csname margininit\endcsname \relax\else\margininit\fi

\newpage

\head {\S1} \endhead  \resetall \sectno=1
\bigskip

The aim of this paper is to prove the following theorem. 
\proclaim{\stag{1.0} Theorem}   If $\kappa$ is a strongly compact cardinal,
$\lambda > \kappa$ is regular and $\zeta,\theta < \kappa$ \ub{then} the partition relation 
$\left(2^{<\lambda}\right)^+ \rightarrow (\lambda + \zeta)^2_\theta$ holds.
\endproclaim
\bigskip

We notice that our argument is valid in the case $\kappa = \omega$. 
As for the history of the problem we notice that  first Hajnal proved 
in an unpublished work, that $\left(2^\omega\right)^+ \rightarrow
(\omega_1 + n)^2_2$ holds for every $n < \omega$. 
Then I showed in \cite{Sh:26} that for $\kappa > \omega$ regular and
$|\alpha|^+ < \kappa$, the relation $\left(2^{< \kappa}\right)^+
\rightarrow (\kappa + \alpha)^2_2$ is true. 
More recently Baumgartner, Hajnal, and Todor\u cevi\'c in
\cite{BHT93} extended this to the case when the number of colors is arbitrary finite. 
\bigskip

\demo{Notation}  If $S$ is a set, $\kappa$ a cardinal then 
$[S]^\kappa = \{a \subseteq S:|a|= \kappa\},[S]^{< \kappa} =
\{a \subseteq S:|a| < \kappa\}$. 
If $D$ is some filter over a set $S$ then $X\in D^+$ denotes that 
$S\setminus X\notin D$. 
If $\kappa < \mu$  are regular cardinals then 
$S^\mu_\kappa = \{\alpha < \mu:\text{cf}(\alpha) =\kappa\}$, a stationary set. 
The notation $A = \{x_\alpha:\alpha < \gamma\}_<$, etc., means that
$A$ is enumerated increasingly. 
\enddemo
\bigskip

\proclaim{\stag{1.1} Lemma}  Assume $\mu = \mu^\theta$.  Assume that
$D$ is a normal filter on $\mu^+$ and $A^* \in D^+,\delta \in A^*
\Rightarrow \text{ cf}(\delta) \ge \lambda$, and $F'$ is a function
with domain $[A]^2$ and range of cardinality $\theta$.
There are a normal filter $D_0$ on $\mu^+$ extending $D,
A_0 \in D_0$ with $A_0 \subseteq A$ and $C_0 \subseteq 
\text{ Rang}(F')$ satisfying Rang$(F' \restriction [A_0]^2) 
= C_0$ such that: if $X \in D^+_0$ \ub{then}
Rang$(F' \restriction [X]^2) \supseteq C_0$.
\endproclaim
\bigskip

\demo{Proof}  Assume indirectly that for no stationary set does the statement of the
Lemma hold.
\enddemo
\bigskip

\proclaim{\stag{1.1a} Claim}   Assume $S^* \subseteq \mu^+$ belongs to
$D^+$ and $\delta \in S_0 \Rightarrow \text{ cf}(\delta_0) >
\theta$.
There is a stationary set $A \subseteq S^*$ and some $C \subseteq
\theta$ such that Rang$(F' \restriction [A]^2)=C$ and: if $f:A \rightarrow
\mu^+$ is a regressive function \ub{then} for some $\alpha < \mu^+$
Rang$(F' \restriction [f^{-1}(\alpha)]^2)=C$ and $f^{-1}(\alpha)$ is a stationary
subset of $\mu^+$ holds.
\endproclaim
\bigskip

\demo{Proof}   Assume that no such sets $A,C$ exist. 
We build a tree $T$ as follows. 
Every node $t$ of the tree will be of the form

$$
\align
t &= \left\langle\langle A_\alpha:\alpha \le \varepsilon \rangle,
\langle f_\alpha:\alpha < \varepsilon \rangle,\langle i_\alpha:\alpha < \varepsilon \rangle
\right\rangle \\
  &= \left< \langle A^t_\alpha:\alpha \le \varepsilon \rangle,\langle
f^t_\alpha:\alpha < \varepsilon \rangle,\langle i^t_\alpha:\alpha <
\varepsilon \rangle \right>
\endalign
$$
\mn
for some ordinal $\varepsilon = \varepsilon(t)$ where $\langle A_\alpha:\alpha \le
\varepsilon \rangle$ is a decreasing, continuous sequence of subsets of $\mu^+$, 
for every $\alpha < \varepsilon,f_\alpha$ is a regressive function on $A_\alpha$, 
and $\langle i_\alpha:\alpha < \varepsilon \rangle$ is a sequence of distinct 
elements of $\theta$. 
It will always be true that if $t<t'$ then each of the three sequences of $t'$ 
extend the corresponding one of $t$.  

To start, we make the node $t$ with $\varepsilon(t)=0,A_0=\mu^+$ the root of the tree. 

At limit levels we extend (the obvious way) all cofinal branches to a node. 

If we are given an element $t= \left\langle \langle A_\alpha:\alpha
\le \varepsilon \rangle,\langle f_\alpha:\alpha < \varepsilon \rangle,
\langle i_\alpha:\alpha < \varepsilon \rangle \right\rangle$ of the
tree and the set $A_\varepsilon$ is $= \emptyset$ mod $D$ then we 
leave $t$ as a terminal node. 
Otherwise, let $C = \text{Rang}\left(F' \restriction 
[A_\varepsilon]^2\right)$ and notice that by 
hypothesis, $A_\varepsilon,C$ toward contradiction, the pair
$A_\varepsilon,C$ cannot be as required in the Claim. 
There is, therefore, a regressive function $f = f^t_t$ with domain
$A_\varepsilon$, such that for every $x < \mu^+$ the set 
Rang$\left(F' \restriction [f^{-1}(x)]^2\right)$ is a $= \emptyset$
mod $D$ subset of $C$ or $f^{-1}(x)$ is a non-stationary subset of $\mu^+$. 
We make the immediate extensions of $t$ the sequences of the form a
$t_x = \left\langle \langle A_\alpha:\alpha \le \varepsilon + 1 \rangle,
\langle f_\alpha:\alpha < \varepsilon + 1 \rangle,\langle i_\alpha:\alpha <
\varepsilon + 1 \rangle \right\rangle$ where $A_{\varepsilon +1} =
f^{-1}(x),f_\alpha= f^t$ and $i_\varepsilon \in C$ is some colour value not in 
the range of $F \restriction [A_\varepsilon]^2$.

Having constructed the tree observe that every element 
$x < \mu^+$ is covered (uniquely) by $A^{t(x)}_{\varepsilon(x)}$ in 
some  terminal node $t(x)$. 
Also, $\varepsilon(x) < \theta^+$ holds by the selection of the $i_\beta$'s. 
For some stationary set $S \subseteq S_0$ of ordinals $x < \mu^+$ the value of
$\varepsilon(x)$ is the same, say $\varepsilon$. 
For $x \in S$ we let $g_\alpha(x) = f^t_\alpha(x)$ where $f^t_\alpha$ is the 
$\alpha$-th regressive function in the node $t$. 
Again, by $\mu^\theta = \mu \and (\forall \alpha \in
S)\text{cf}(\alpha) > \theta$ we have that $(\forall x \in S') \dsize
\wedge_{\alpha < \varepsilon} g_\alpha(x) = \beta_\alpha$ holds for a 
subset $S'\subseteq S$ for $D^+$. 
But then we get that the stationary set $S'$ 
satisfies $x,y \in S' \Rightarrow
(A^{t(x)}_\alpha,f^{t(x)}_\alpha,i^{t(*)}_\alpha) =
(A^{t(y)}_\alpha,f^{t(y)}_\alpha,i^{t(y)}_\alpha)$; we can prove this
by induction on $\alpha$ and we can conclude that $x,y \in S'
\Rightarrow t(x) = t(y)$, so $S' \subseteq A^t_{\varepsilon(t)}$ for
some terminal node $t$, so this latter set is in $D^+$, a
contradiction.  \nl
${{}}$  \hfill$\square_{\scite{1.1a}}$\margincite{1.1a}
\enddemo
\bigskip

\demo{Continuation of the proof of Lemma \scite{1.2}}  Define the ideal $I$ as follows. 
For $X \subseteq \mu^+$ we let $X \in I$ \ub{iff} there are a club $E$ of
$\mu^+$ and a regressive $f:X \cap A \rightarrow \mu^+$ such that every 
Rang$\left(F' \restriction [f^{-1}(\alpha)]^2\right)$ is a proper subset of $C$ or
$f^{-1}(\alpha)$ is a $= \emptyset$ mod $D$ subset of $\mu^+$. 
\enddemo
\bigskip

\proclaim{\stag{1.1b} Claim}   $I$ is a normal ideal on $\mu^+$.
\endproclaim
\bigskip

\demo{Proof}  Straightforward. 

Now apply Claim \scite{1.1a} for the $S^* = A^*$ to get $(C,A)$ and
then define $I$ as above: set $D_0$ to be the dual filter 
of $I$, let $A_0 = A$ and let $C_0 = C$; by \scite{1.1b} we are done.
\hfill$\square_{\scite{1.1}}$\margincite{1.1}
\enddemo
\bigskip

\remark{\stag{1.1c} Remark}   1) If Lemma \scite{1.1} holds 
for some $D_0,A_0,C_0$ then it holds 
for $D_1,A_1,C_0$ when the normal filter $D_1$ extends $D_0$, and 
$A_1\in D_1$ satisfies $A_1 \subseteq A_0$. \nl
2)  If $D_0,A_0,C_0$ satisfy Lemma \scite{1.2}, and 
$X \in D^+_0$ \ub{then} $X$ contains a homogeneous 
set of order type $\lambda + 1$ of color $\xi$ for every  $\xi \in
C_0$. \nl
3) This is closely related to the proof in \cite{Sh:26}.
\endremark
\bigskip

\demo{Proof}  Let $\mu=2^{< \lambda}$, and $F:[\mu^+]^2 \rightarrow
\theta$ be a colouring we apply \scite{1.1} for $S =
S^{\mu^+}_\lambda$, $(\theta = \theta,\mu = \mu)$ and $D$ the club filter. 
\nl
We fix $A_0,D_0,C_0$ as in \scite{1.1}.
\enddemo
\bigskip

\proclaim{\stag{1.2} Lemma}    Almost every $\delta \in A_0$; (i.e. for
all but non-stationarily many) satisfies the following: 
if $s \in [A^*_0 \cap \delta]^{<\lambda}$ and 
$\{z_\alpha:\alpha < \gamma\}_< \subseteq A^* \cap
(\delta,\mu^+)$ and $\gamma < \kappa$ \ub{then} there is 
$\{y_\alpha:\alpha < \gamma\}_< \subseteq A_0 \cap \left(
\sup(s),\delta \right)$ such that:
\mr
\item "{$(a)$}"   $F(x,y_\alpha) = F(x,z_\alpha) \quad$ (for $x\in
s,\alpha < \gamma$);
\sn
\item "{$(b)$}"  $F(y_\alpha,y_\beta) = F(z_\alpha,z_\beta) \quad$ 
(for $\alpha < \beta < \gamma$). 
\endroster
\endproclaim
\bigskip

\demo{Proof}  By simple reflection.
\enddemo
\bigskip

\proclaim{\stag{1.3} Lemma}   There is $A'_0 \subseteq A_0,A'_0 \in D_0$
such that: if $\delta \in A'_0,s \in[\delta]^{<\lambda}$ and $\varepsilon \in C_0$, 
\ub{then} there exists a $\delta_1 \in A_0,\delta < \delta_1$ such that
\mr
\item "{$(a)$}"  $F(x,\delta) = F(x,\delta_1) \quad$ (for $x\in s$);
\sn
\item "{$(b)$}"  $F(\delta,\delta_1) = \varepsilon$.
\endroster
\endproclaim
\bigskip

\demo{Proof}   Otherwise, there is some $X \subseteq A_0,X \in D^+_0$ such that 
for every $\delta \in X$ there are $s(\delta) \in [\delta]^{<\lambda}$ and 
$\xi(\delta) \in C_0$ such that there is no $\delta_1 > \delta$ satisfying (a) and (b). 
By normality we can assume that $s(\delta) = s$ and $\xi(\delta) = \xi$ holds 
for $\delta \in X$.  By Lemma \scite{1.2} there must exist $\delta < \delta_1$
in $X$ with $F(\delta,\delta_1) = \xi$ and this is a
contradiction. \nl
${{}}$  \hfill$\square_{\scite{1.3}}$
\enddemo
\bn
Now let $A'_0$ satisfy Lemmas \scite{1.1} and \scite{1.3} 
and pick some $\delta_1 \in A'_0$ and then let 
$T = A'_0 \backslash (\delta_1 +1)$. 
\proclaim{\stag{1.4} Lemma}    There exists a function $G:T \times T
\rightarrow C_0$ such that: if $s \in [\delta_1]^{<\lambda},\gamma <
\kappa$, and $Z = \{z_\alpha:\alpha < \gamma\}_< \subseteq T$
\ub{then} there is $\{y_\alpha:\alpha < \gamma\}_< \subseteq \left(\sup(s),\delta_1
\right)$ such that
\mr
\item "{$(a)$}"   $F(x,y_\alpha) = F(x,z_\alpha) \quad$ (for $x \in
s,\alpha < \gamma$);
\sn
\item "{$(b)$}"   $F(y_\alpha,y_\beta) = F(z_\alpha,z_\beta) \quad$ 
(for $\alpha < \beta < \gamma$);
\sn
\item "{$(c)$}"   $F(y_\alpha,z_\beta) = G(z_\alpha,z_\beta) \quad$ 
(for $\alpha,\beta < \gamma$).
\endroster
\endproclaim
\bigskip

\demo{Proof}  As $\kappa$ is strongly compact, it suffices to show that for every 
$Z \in [T]^{< \kappa}$ there exists a function $G:Z \times Z
\rightarrow \theta$ as required. 
Clauses $(a)$ and $(b)$ are obvious by Lemma \scite{1.2}, and it is clear that, if we fix 
$Z$, then for every  $s \in [\delta_1]^{<\lambda}$ there is an appropriate 
$G:Z \times Z \rightarrow \theta$. 
We show that there is some $G:Z \times Z \rightarrow \theta$ that works for every 
$s$. 
Assume otherwise, that is, for every $G:Z \times Z \rightarrow \theta$ there is 
some $s_G \in [\delta_1]^{<\lambda}$ such that $G$ is not appropriate for $s_G$. 
Notice that the number of these functions $G$ is less than $\kappa$. 
Then no $G$ could be right for $s = \cup\{s_G:G$ a function from $Z
\times Z$ to $\theta\} \in [\delta_1]^{<\lambda}$, a contradiction.
\hfill$\square_{\scite{1.4}}$\margincite{1.4}
\enddemo
\bigskip

\demo{Continuation of the proof of Theorem \scite{1.0}}
We now apply Lemma \scite{1.1} to the coloring $\bar G(x,y) = 
\left\langle F(x,y),G(y,x) \right\rangle$ for $x < y$ in $T$ and $0$
otherwise, and the filter $D_0$ and the set $T$ and get the
normal filter $D_1 \supseteq D_0$, the set $A_1 \subseteq T,A_1 \in
D_1$ and the colour set $C_1 \subseteq \theta \times \theta$. 
Notice that actually $C_1 \subseteq C_0 \times C_0$. 
We can also apply Lemmas \scite{1.2} and \scite{1.3} and get some set
$A'_1 \subseteq A_1$. 
\enddemo
\bigskip

\proclaim{\stag{1.5} Lemma}    There is a set $a \in [A'_1]^{< \kappa}$
such that for every decomposition $a = \cup\{a_\xi:\xi < \theta\}$
there is some $\xi <\theta$ such that for every $\varepsilon \in C_1$ there is 
an $\varepsilon$-homogeneous subset for the colouring $\bar G$ of order type
$\zeta$ in $a_\xi$.
\endproclaim
\bigskip

\demo{Proof}  This follows from the strong compactness of $\kappa$ 
as $A'_1$ itself has this partition
property. \hfill$\square_{\scite{1.5}}$\margincite{1.5}
\enddemo
\bn
Fix a set $a$ as in \scite{1.5}.

We now describe the construction of the required homogeneous subset. 
Let $\delta_2 \in A'_1$ be some element with $\delta_2 > \sup(a)$. 
For $\xi < \theta$ let $a_\xi$ be the following set: 

$$
a_\xi = \{x \in a:G(\delta_2,x) = \xi\}.
$$
\mn
By Lemma \scite{1.5}, there is some $\xi <\theta$ for 
which the statement as above is true and necessarily (as $s \cup
\{\delta_2\} \subseteq A'_1 \subseteq A_0$) we have $\xi \in C_0$.        
We notice that $\langle \xi,\xi' \rangle \in C_1$ holds for an appropriate 
$\xi'\in C_0$. 
Select some $b \subseteq a_\xi$, tp$(b)= \zeta$ such that $F$ is constantly  
$\xi$ on $b$. 
This set $b$ will be the $\zeta$ part of our homogeneous set of ordinals
$\lambda + \zeta$, so we will have to construct a set of order type $\lambda$ below 
$b$. 
By induction on $\alpha$ we will choose $x_\alpha$ such that the set 
$\{x_\alpha:\alpha < \lambda\}_< \subseteq \delta_1$ satisfies the following conditions: 
\smallskip
$F(x_\beta,x_\alpha) = \xi$ (for $\beta < \alpha$), 
\smallskip
$F(x_\alpha,b \cup \{\delta_2\}) = \xi$, i.e. $F(x_\alpha,y) = \xi$
when $y \in b \cup \{\delta_2\}$.

Assume that we have reached step $\alpha$, that is, we are given the 
ordinals $\{x_\beta:\beta < \alpha\}_<$.  
Applying Lemma \scite{1.3} we get that there exists some $\gamma > \delta_2$ 
such that $\gamma \in A_1$ moreover 
\smallskip
$F(x_\beta,\gamma) = \xi$ (for $\beta < \alpha$), 
\smallskip
$G \left(\gamma,b \cup \{\delta_2\} \right) = \xi$, i.e. $G(\gamma,y)
= \xi$ when $y \in b \cup \{\delta_2\}$. 
\mn
By the definition of $G$ this implies that there is some $x_\alpha$ as required. 
\hfill$\square_{\scite{1.0}}$\margincite{1.0}

\newpage
    
REFERENCES.  
\bibliographystyle{lit-plain}
\bibliography{lista,listb,listx,listf,liste}

\enddocument

\bye